\documentclass[12pt,a4paper]{article}
\usepackage{amsfonts}
\usepackage{amssymb}
\usepackage{mathrsfs}
\usepackage{amsmath}
\usepackage{amsthm}
\usepackage{enumerate}
\usepackage{cite}
\usepackage[all]{xy}
\usepackage{extarrows}
\usepackage{indentfirst}
\allowdisplaybreaks
\newtheorem{defn}{Definition}[section]
\newtheorem{thm}[defn]{Theorem}
\newtheorem{lem}[defn]{Lemma}
\newtheorem{prop}[defn]{Proposition}
\newtheorem{cor}[defn]{Corollary}
\newtheorem{eg}[defn]{Example}
\newtheorem{re}[defn]{Remark}

\newcommand\relphantom[1]{\mathrel{\phantom{#1}}}
\newcommand{\bdefn}{\begin{defn}}
\newcommand{\edefn}{\end{defn}}
\newcommand{\bthm}{\begin{thm}}
\newcommand{\ethm}{\end{thm}}
\newcommand{\blem}{\begin{lem}}
\newcommand{\elem}{\end{lem}}
\newcommand{\bprop}{\begin{prop}}
\newcommand{\eprop}{\end{prop}}
\newcommand{\bcor}{\begin{cor}}
\newcommand{\ecor}{\end{cor}}
\newcommand{\beg}{\begin{eg}}
\newcommand{\eeg}{\end{eg}}
\newcommand{\bre}{\begin{re}}
\newcommand{\ere}{\end{re}}
\newcommand{\bpf}{\begin{proof}}
\newcommand{\epf}{\end{proof}}
\newcommand{\benu}{\begin{enumerate}}
\newcommand{\eenu}{\end{enumerate}}
\newcommand{\bc}{\begin{center}}
\newcommand{\ec}{\end{center}}
\newcommand{\bea}{\begin{eqnarray}}
\newcommand{\eea}{\end{eqnarray}}
\newcommand{\Bea}{\begin{eqnarray*}}
\newcommand{\Eea}{\end{eqnarray*}}
\newcommand{\beq}{\begin{equation}}
\newcommand{\eeq}{\end{equation}}
\newcommand{\Beq}{\begin{equation*}}
\newcommand{\Eeq}{\end{equation*}}
\newcommand{\bspl}{\begin{split}}
\newcommand{\espl}{\end{split}}

\begin{document}

\title{\textbf{On  split  Regular Hom-Lie color algebras}
\author{ Yan Cao$^{1,2},$  Liangyun Chen$^{1}$
 \date{{\small {$^1$ School of Mathematics and Statistics, Northeast Normal
 University,\\
Changchun 130024, China}\\{\small {$^2$  Department of Basic
 Education,
 Harbin University of
Science and Technology,\\ Rongcheng Campus,  Rongcheng 264300,
China}}}}}} \maketitle
\date{}

\begin{abstract}

We introduce the class of split regular Hom-Lie color algebras as the natural generalization of  split Lie color algebras. By developing techniques of connections of roots for this kind of algebras,  we show that  such a split regular Hom-Lie color algebra  $L$ is of the form
$L = U + \sum\limits_{[j] \in \Lambda/\sim}I_{[j]}$ with $U$ a subspace of the abelian graded subalgebra $H$ and any $I_{[j]}$, a well described
ideal of $L$, satisfying $[I_{[j]}, I_{[k]}] = 0$ if $[j]\neq [k]$. Under certain conditions, in the case of $L$ being of maximal length,  the simplicity of the algebra is characterized. \\

\end{abstract}
\renewcommand{\thefootnote}{\fnsymbol{footnote}}
\footnote[0]{ Corresponding author(L. Chen): chenly640@nenu.edu.cn.}
\footnote[0]{Supported by  NNSF of China (Nos. 11171055 and
11471090),  Scientific
Research Fund of Heilongjiang Provincial Education Department
 (No. 12541184). }

\section{Introduction}
The notion of Hom-Lie algebras was introduced by Hartwig, Larsson and Silvestrov to describe the $q$-deformation of the Witt and the Virasoro algebras \cite {1}. Since then, many authors have studied Hom-type algebras  \cite{2,3,4,5,6,7}. The notion of Lie color algebras was introduced as generalized Lie algebras in 1960 by Ree  \cite{L2}. So far, many results of this kind of algebras have been considered in the frameworks of enveloping
algebras, representations and related problems \cite{AAO, AAO2, AO}. In particular, Yuan introduced the notion of  Hom-Lie color algebras in \cite {8}, which can be
viewed as an extension of Hom-Lie (super)algebras to $\Gamma$-graded algebras, where $\Gamma$ is any
abelian group.

As is well-known, the class of the split algebras is specially related to addition quantum
numbers, graded contractions and deformations. For instance, for a physical system which displays
a symmetry of $L$, it is interesting to know in detail the structure of the split decomposition
because its roots can be seen as certain eigenvalues which are the additive quantum numbers
characterizing the state of such system. Determining the structure of split algebras will become more and more meaningful in the area of research in mathematical physics. Recently, in \cite{BL52, BL52567, BL5234, BL528}, the structure of arbitrary split Lie algebras, arbitrary split
  Lie color algebras, arbitrary split  Lie triple systems and arbitrary split regular Hom-Lie algebras have been determined by the techniques of connections of roots.
 The purpose of this paper is to consider the structure of split regular Hom-Lie color algebras by the techniques of connections of roots based on some work in \cite{BL5234, BL52567}.

Throughout this paper, split regular Hom-Lie color algebras $L$ are considered of
arbitrary dimension and over an arbitrary base field $\mathbb{K}$.  This paper
is organized as follows. In section 2, we establish the preliminaries on
split regular Hom-Lie color algebras theory.  In section 3, we show that such
an arbitrary  regular Hom-Lie color algebra  $L$ with a symmetric root system is
of the form  $L=U+\sum_{[j]\in \Lambda/\sim} I_{[j]}$ with $U$ a
subspace of the abelian  subalgebra $H$  and any $I_{[j]}$ a well described ideal of
$L$,  satisfying $[I_{[j]},I_{[k]}]
=0$ if $[j]\neq [k]$. In section 4, we show that under certain conditions, in
the case of $L$ being of maximal length,  the simplicity
of the algebra is characterized.

\section{Preliminaries}
First we recall the definitions of  Lie color algebras and Hom-Lie color algebras. The following definition is well-known from the theory of graded algebra.
\bdefn{\rm\cite{8}} Let $\Gamma$ be an abelian group. A bi-character on $\Gamma$ is a map $\varepsilon:\Gamma \times \Gamma \rightarrow \mathbb{K}\setminus \{0\}$ satisfying

$\rm 1$. $\varepsilon(\alpha,\beta)\varepsilon(\beta,\alpha)=1,$

$\rm 2$. $\varepsilon(\alpha,\beta+\gamma)=\varepsilon(\alpha,\beta)\varepsilon(\alpha,\gamma),$

$\rm 3$. $\varepsilon(\alpha+\beta,\gamma)=\varepsilon(\alpha,\gamma)\varepsilon(\beta,\gamma),$

\noindent for all $\alpha, \beta, \gamma \in \Gamma$.
\edefn

 It is clear that $\varepsilon(\alpha,0)=\varepsilon(0,\alpha)=1$ for any $\alpha \in \Gamma$, where 0 denotes the identity element of $\Gamma$.

\bdefn{\rm\cite{BL5234}} Let $L=\oplus_{g \in \Gamma}L_{g}$ be a $\Gamma$-graded $\mathbb{K}$-vector space. For a nonzero homogeneous element $v \in L$,
denote by $\bar{v}$ the unique group element in $\Gamma$ such that $v \in L_{\bar{v}}$, which will be called the homogeneous degree of $v$. We shall say that $L$ is a \textbf{Lie color algebra} if it is endowed with a $\mathbb{K}$-bilinear map
$[\cdot,\cdot]:L\times L \rightarrow L$ satisfying

$[v,w]=-\varepsilon(\bar{v},\bar{w})[w,v],$ $(\mathrm{Skew}$- $\mathrm{symmetry})$

$[v,[w,t]]=[[v,w],t]+\varepsilon(\bar{v},\bar{w})[w,[v,t]],$ $(\mathrm{Jacobi}$ $\mathrm{identity})$

\noindent for all homogeneous elements $v,w,t \in L$.
\edefn

 Lie superalgebras are examples of Lie color algebras  with $\Gamma=\mathbb{Z}_{2}$ and $\varepsilon(i,j)=(-1)^{ij}$, for any $i,j \in \mathbb{Z}_{2}$. We also note that $L_{0}$ is a Lie algebra.

\bdefn{\rm\cite{8}}  A \textbf{Hom-Lie color algebra} is a quadruple $(L, [\cdot,\cdot],\phi,\varepsilon)$ consisting of a $\Gamma$-graded space $L$, an even bilinear mapping  $[\cdot,\cdot]: L\times L\rightarrow L$, a homomorphism $\phi$ and a bi-character $\varepsilon$ on $\Gamma$ satisfying

$[x,y]=-\varepsilon(\bar{x},\bar{y})[y,x],$

$\varepsilon(\bar{z},\bar{x})[\phi(x),[y,z]]+\varepsilon(\bar{x},\bar{y})[\phi(y),[z,x]]+\varepsilon(\bar{y},\bar{z})[\phi(z),[x,y]]=0,$

\noindent for all homogeneous elements $x,y,z \in L$,  $\bar{x},\bar{y},\bar{z}$ denote the homogeneous degree of $x,y,z$.
 When $\phi$ furthermore is an algebra automorphism it is said that $L$ is a \textbf{regular Hom-Lie color algebra}.
\edefn

Clearly Hom-Lie algebras and Lie color algebras are examples of Hom-Lie color algebras.

Throughout this paper we will consider regular Hom-Lie color algebras $L$ being of arbitrary dimension and arbitrary base field $\mathbb{K}$.  $\mathbb{N}$ denotes  the set of all non-negative integers and  $\mathbb{Z}$ denotes the set of all integers.

For any $x\in L$, we consider the adjoint mapping $\mathrm{ad}_{x}:L\rightarrow L$ defined by $\mathrm{ad}_{x}(z)=[x,z].$
The usual regularity concepts will be understood in the graded sense. For instance,
A subalgebra $A$ of $L$ is a graded subspace $I=\oplus_{g \in \Gamma}I_{g}$ of $L$ such that $[A,A] \subset A$ and $\phi(A) =
A$. A graded subspace $I=\oplus_{g \in \Gamma}I_{g}$ of $L$ is called an ideal if $[I,L] \subset I$ and $\phi(I) = I$. A Hom-Lie color algebra $L$ will be called simple if
$[L, L] \neq 0$ and its only (graded) ideals are {0} and $L$.

We introduce the concept of split regular Hom-Lie  color algebra in an analogous way. We begin by considering a
maximal abelian graded subalgebra $H =\oplus_{g \in \Gamma}H_{g}$  among the abelian graded subalgebras of $L$. Observe
that $H$ is necessarily a maximal abelian subalgebra of $L$ as the following lemma shows.

 \blem \label{1.1}
 Let $H =\oplus_{g \in \Gamma}H_{g}$ be a maximal abelian graded subalgebra of a Hom-Lie color algebra $L$. Then $H$
is a maximal abelian subalgebra of $L$.

 \elem
\bpf
Consider an abelian subalgebra $K$ of $L$ such that $H \subset K$. For any $x \in K$ we have $[x, H_{g}] = 0$
for each $g \in \Gamma$, and so by writing $x =\sum_{i=1}^{n} x_{g_{i}}$ with $x_{g_{i}} \in L_{g_{i}}$ for $i = 1, \cdots, n$, being $g_{i} \in \Gamma$
and $g_{i}\neq g_{j}$ if $i \neq j$, we get by the grading $[x_{g_{i}} , H_{g}] = 0$. Hence, for any $g_{i}$, $i = 1, \cdots n$, we have
$(H_{g_{i}}
+ \mathbb{K}x_{g_{i}} ) \oplus(\oplus_{g \in \Gamma\setminus \{g_{i}\}}H_{g})$ is an abelian graded subalgebra of $L$ containing $H$ and so $x_{g_{i}} \in H_{g_{i}}$.
 From here we get $x \in H$ and then $K = H$.
\epf

Let us introduce the class of split algebras in the framework of regular Hom-Lie color
algebras. Denote by $H =\oplus_{g \in \Gamma}H_{g}$ a maximal  abelian (graded) subalgebra, (MAGSA), of a regular Hom-Lie color algebra $L$. For a linear
functional
$$\alpha:H_{0}\rightarrow \mathbb{K},$$
we define the root space of $L$ $($with respect to $H$$)$ associated to $\alpha$ as the subspace
$$L_{\alpha} = \{v_{\alpha}\in L: [h_{0}, v_{\alpha} ] = \alpha(h_{0})\phi(v_{\alpha}) \ for \ any \ h_{0} \in H_{0}\}.$$
The elements $\alpha:H_{0}\rightarrow \mathbb{K}$ satisfying $L_{\alpha}\neq 0$ are called roots of $L$ with respect to
$H$. We denote $\Lambda:= \{\alpha \in H_{0}^{\ast}
\setminus \{o\}: L_{\alpha}\neq 0\}$.

\bdefn
We say that $L$ is a \textbf{split regular Hom-Lie color algebra}, with respect
to $H$, if
$$L=H\oplus (\oplus_{\alpha \in \Lambda }L_{\alpha}).$$
We also say that $\Lambda$ is the root system of $L$.
\edefn

Note that when $\phi = \mathrm{Id}$, the split Lie color algebras become examples of split
regular Hom-Lie color algebras. Hence, the present paper extends the results in \cite{BL5234}.
For convenience, the mappings $\phi| _{H}$, $\phi|_{H}^{-1}:H \rightarrow H$ will be denoted by $\phi$ and $\phi^{-1}$ respectively.

It is clear that  the root space associated to the zero root $L_{o}$ satisfies $H \subset L_{o}$. Conversely, given any
$v_{o} \in L_{o}$ we can write $$v_{o} =h\oplus(\oplus_{i=1}^{n}v_{\alpha_{i}}),$$
 where $h \in H$ and $v_{\alpha_{i}} \in L_{\alpha_{i}}$
for $i = 1,\cdots, n,$
with $\alpha_{i} \neq \alpha_{j}$ if $i \neq j$. Hence $$0 = [h_{0}, h \oplus(\oplus_{i=1}^{n}v_{\alpha_{i}}) ]=\oplus_{i=1}^{n}\alpha_{i}(h_{0}) \phi(v_{\alpha_{i}}),$$ for any $h_{0}\in H_{0}$. So taking into account the direct character of the sum and that $\alpha_{i}\neq 0$ give us  $v_{\alpha_{i}}=0$ for $i = 1,\cdots, n$. So $v_{o}=h \in H$. Consequently,
\begin{equation}\label{zhang}
H=L_{o}.
\end{equation}
 \blem \label{1.2}
 Let $L=\oplus_{g \in \Gamma}L_{g}$ be a split Hom-Lie color algebra with corresponding root space decomposition $L=H\oplus (\oplus_{\alpha \in \Lambda }L_{\alpha}).$
If we denote by $L_{\alpha,g}=L_{\alpha}\cap L_{g}$, then the following assertions hold.

$\rm 1$. $L_{\alpha}=\oplus_{g \in \Gamma}L_{\alpha,g}$ for any $\alpha \in \Lambda\cup \{o\}$.

$\rm 2$. $H_{g}=L_{o,g}.$ In particular $H_{0}=L_{o,0}$.

$\rm 3$. $L_{0}$ is a split Hom-Lie algebra, respect to $H_{0}$, with root space decomposition $L_{0}=H_{0}\oplus (\oplus_{\alpha \in \Lambda }L_{\alpha,0}).$
\elem

\bpf $\rm 1$. By the $\Gamma$-grading of $L$ we may express any $v_{\alpha}\in L_{\alpha}$, $\alpha \in \Lambda\cup \{o\}$, in the form $v_{\alpha}=v_{\alpha,g_{1}}+\cdots+v_{\alpha,g_{n}}$ with $v_{\alpha,g_{i}} \in L_{g_{i}}$ for distinct $g_{1},\cdots,g_{n} \in \Gamma$. If $h_{0}\in H_{0}$ then $[h_{0},v_{\alpha,g_{i}}]=\alpha(h_{0})\phi(v_{\alpha,g_{i}})$ for $i=1,\cdots,n.$ Hence $L_{\alpha}=\oplus_{g \in \Gamma}(L_{\alpha}\cap L_{g})$ and we can write $L_{\alpha}=\oplus_{g \in \Gamma}L_{\alpha,g}$ for any $\alpha \in \Lambda\cup \{o\}$.

$\rm 2$. Consequence of Eq. (\ref{zhang}) and item 1.

$\rm 3$. We also have $L_{g}=H_{g}\oplus(\oplus_{\alpha \in \Lambda}L_{\alpha,g})$ for any $g\in  \Gamma$. By considering $g=0$ we get $L_{0}=H_{0}\oplus(\oplus_{\alpha \in \Lambda}L_{\alpha,0}).$ Hence, the direct character of the sum and the fact that $\alpha \neq 0$ for any $\alpha \in \Lambda$ give us that $H_{0}$ is a MASA of the  Hom-Lie algebra $L_{0}$. Hence $L_{0}$ is a split Hom-Lie algebra respect to  $H_{0}$.
\epf

 \blem \label{3556000}
For any $\alpha$, $\beta \in \Lambda \cup \{o\}$, the following assertions hold.

$\rm 1$. $\phi(L_{\alpha})\subset L_{\alpha \phi^{-1}}$ and $\phi^{-1}(L_{\alpha})\subset L_{\alpha \phi}$.

$\rm 2$. $[L_{\alpha}, L_{\beta}] \subset L_{\alpha\phi^{-1}+\beta\phi^{-1}}$.
\elem

\bpf $\rm 1$.  For $h_{0} \in H_{0}$ write $h_{0}^{'} = \phi(h_{0})$. Then for all $h_{0} \in H_{0}$ and $v_{\alpha} \in L_{\alpha}$, since
$[h_{0}, v_{\alpha}] = \alpha(h_{0})\phi(v_{\alpha})$, one has
$$[h_{0}^{'}, \phi(v_{\alpha})] = \phi([h_{0}, v_{\alpha}]) = \alpha(h_{0})\phi(\phi(v_{\alpha})) =\alpha\phi^{-1}(h_{0}^{'})\phi(\phi(v_{\alpha})).$$
Therefore we get $\phi(v_{\alpha}) \in L_{\alpha\phi^{-1}}$ and so $\phi(L_{\alpha}) \subset L_{\alpha \phi^{-1}}$. In a similar way, one gets
$\phi^{-1}(L_{\alpha})\subset L_{\alpha \phi}$.

$\rm 2$. For any $h_{0} \in H_{0}$, $v_{\alpha} \in L_{\alpha}$ and $v_{\beta} \in L_{\beta}$, by denoting $h_{0}^{'} = \phi(h_{0})$, by Hom-Jacobi identity, we
have that \begin{align*}
[h_{0}^{'}, [v_{\alpha}, v_{\beta}] ]&= [[h_{0}, v_{\alpha}],\phi(v_{\beta})]+\varepsilon(\bar{h_{0}},\bar{v_{\alpha}})[\phi(v_{\alpha}),[h_{0},v_{\beta}]]\\
&=[\alpha(h_{0})\phi(v_{\alpha}),\phi(v_{\beta})]+\beta(h_{0})[\phi(v_{\alpha}), \phi(v_{\beta})]\\
&= (\alpha+\beta)(h_{0})\phi([v_{\alpha}, v_{\beta}])\\
&= (\alpha+\beta)\phi^{-1}(h_{0}^{'})\phi([v_{\alpha}, v_{\beta}]).
\end{align*}
Therefore we get $[v_{\alpha}, v_{\beta}] \in L_{\alpha\phi^{-1}+\beta\phi^{-1}}$ and  so  $[L_{\alpha}, L_{\beta}] \subset L_{\alpha\phi^{-1}+\beta\phi^{-1}}$.
\epf

 \blem \label{3556000asd}
If $\alpha \in \Lambda$ then $\alpha\phi^{-z} \in \Lambda$ for any $z\in \mathbb{Z}.$
\elem

\bpf  It is a consequence of Lemma \ref{3556000}-1.
\epf

\bdefn
 A root system $\Lambda$  of a split Hom-Lie color algebra  is called \textbf{symmetric} if it satisfies that
$\alpha \in \Lambda$ implies $-\alpha \in \Lambda$.
\edefn

\section{Decompositions}
In the following, $L$ denotes a split regular  Hom-Lie color algebra with a symmetric root system $\Lambda$ and $L=
H \oplus(\oplus _{\alpha \in \Lambda}L_{\alpha})$ the corresponding root decomposition. Given a linear functional  $\alpha:H\rightarrow \mathbb{K}$, we denote by $-\alpha:H\rightarrow \mathbb{K}$ the element in $H^{\ast}$ defined by $(-\alpha)(h):=-\alpha(h)$ for all $h\in H$. We begin by developing the techniques of  connections of roots
 in this section.

\bdefn\label{egzmmm}
Let $\alpha$ and $\beta$ be two nonzero roots. We shall say that $\alpha$ is
\textbf{connected} to $\beta$  if there exists $\alpha_{1},\cdots,\alpha_{k}\in \Lambda$  such that

 $\rm 1$. If $k=1$, then
 $\alpha_{1}\in \{a\phi^{-n}:  n\in \mathbb{N}\}\cap \{\pm \beta \phi^{-m}:  m\in \mathbb{N}\}.$
If $k\geq 2$, then
  $\alpha_{1}\in \{a\phi^{-n}:  n\in \mathbb{N}\}$.

$\rm 2$. $\alpha_{1}\phi^{-1}+\alpha_{2}\phi^{-1} \in \Lambda$,

\  \quad  $\alpha_{1}\phi^{-2}+\alpha_{2}\phi^{-2}+\alpha_{3}\phi^{-1} \in \Lambda$,

\  \quad  $\alpha_{1}\phi^{-3}+\alpha_{2}\phi^{-3}+\alpha_{3}\phi^{-2}+\alpha_{4}\phi^{-1} \in \Lambda$,

\  \quad $\cdots \cdots \cdots$

\  \quad $\alpha_{1}\phi^{-i}+\alpha_{2}\phi^{-i}+\alpha_{3}\phi^{-i+1}+\cdots+\alpha_{i+1}\phi^{-1}\in \Lambda$,

\  \quad $\cdots \cdots \cdots$

\  \quad $\alpha_{1}\phi^{-k+2}+\alpha_{2}\phi^{-k+2}+\alpha_{3}\phi^{-k+3}+\cdots+\alpha_{i}\phi^{-k+i}+\cdots+\alpha_{k-1}\phi^{-1}\in \Lambda$.

$\rm 3$. $\alpha_{1}\phi^{-k+1}+\alpha_{2}\phi^{-k+1}+\alpha_{3}\phi^{-k+2}+\cdots+\alpha_{i}\phi^{-k+i-1}+\cdots+\alpha_{k}\phi^{-1}\in \{\pm \beta \phi^{-m}: m \in \mathbb{N}\}$.

\noindent We shall also say
that $\{\alpha_{1},\cdots,\alpha_{k}\}$ is a connection from $\alpha$ to $\beta$.
\edefn

Observe that the case $k =1$ in Definition \ref{egzmmm}  is equivalent to the fact
$\beta=\epsilon \alpha\phi^{z}$ for some $z \in \mathbb{Z}$ and $\epsilon \in \{\pm 1\}$.

 \blem \label{3556jjsd}
For any $\alpha \in \Lambda$, we have that $\alpha\phi^{z_{1}}$ is connected to  $\alpha\phi^{z_{2}}$ for every
$z_{1}, z_{2} \in \mathbb{Z}.$ We also have that $\alpha \phi^{z_{1}}$ is connected to $-\alpha \phi^{z_{2}}$ in case  $-\alpha \phi^{z_{2}}\in \Lambda$.

\elem

\bpf
This can be proved completely analogously to  \cite[Lemmas 2.2]{BL52567}.
\epf

\blem \label{3556jjsd77777}
Let $\{\alpha_{1},\cdots,\alpha_{k}\}$ be a connection from   $\alpha$ to  $\beta$. Then the following
assertions hold.

$\rm 1$. Suppose $\alpha_{1}=\alpha \phi^{-n}, n\in \mathbb{N}$. Then for any $r\in \mathbb{N}$  such that $r\geq n$, there
exists a connection $\{\overline{\alpha}_{1},\cdots,\overline{\alpha}_{k}\}$ from $\alpha$ to $\beta$ such that $\overline{\alpha}_{1}=\alpha \phi^{-r}$.

$\rm 2$. Suppose that $\alpha_{1}=\epsilon\beta\phi^{-m}$ in case $k =1$ or
$$\alpha_{1}\phi^{-k+1}+\alpha_{2}\phi^{-k+1}+\alpha_{3}\phi^{-k+2}+\cdots+\alpha_{k}\phi^{-1}= \epsilon\beta \phi^{-m}$$
in case $k \geq 2$, with $m \in \mathbb{N}$ and $\epsilon \in \{\pm 1\}$. Then for any $r\in \mathbb{N}$ such
that $r\geq m$, there exists a connection  $\{\overline{\alpha}_{1},\cdots,\overline{\alpha}_{k}\}$ from $\alpha$ to $\beta$ such that
 $\overline{\alpha}_{1}=\epsilon\beta \phi^{-r}$ in case $k = 1$ or
$$\overline{\alpha}_{1}\phi^{-k+1}+\overline{\alpha}_{2}\phi^{-k+1}+\overline{\alpha}_{3}\phi^{-k+2}+\cdots+\overline{\alpha}_{k}\phi^{-1}= \epsilon\beta \phi^{-r}$$
in case $k \geq 2$.
\elem

\bpf
This can be proved completely analogously to  \cite[Lemmas 2.3]{BL52567}.
\epf

\bprop\label{alp93556jjsd77777}
The relation $\sim$ in $\Lambda$, defined by $\alpha \sim \beta$ if and only if $\alpha$ is connected to $\beta$, is of equivalence.

\eprop

\bpf
This can be proved completely analogously to  \cite[Proposition 2.4]{BL52567}.
\epf

 For any $\alpha \in \Lambda$, we denote by
 $$\Lambda_{\alpha}:=\{\beta \in \Lambda : \beta\sim \alpha\}.$$

\noindent Clearly if $\beta \in \Lambda_{\alpha}$ then $-\beta \in \Lambda_{\alpha}$ and, by Proposition  \ref{alp93556jjsd77777}, if $\gamma \not \in \Lambda_{\alpha}$ then $\Lambda_{\alpha}\cap \Lambda_{\gamma}=\emptyset$.

 Our next goal is to associate an adequate ideal $L_{\Lambda_{\alpha}}$  of $L$ to any $\Lambda_{\alpha}$. For  $\Lambda_{\alpha}$, $\alpha \in \Lambda$, we define
$$H_{\Lambda_{\alpha}}:=\mathrm{span_{\mathbb{K}}}\{[L_{\beta}, L_{-\beta}]: \beta \in \Lambda_{\alpha}\}.$$
Then $H_{\Lambda_{\alpha}}$ is the direct sum of
$$\sum_{\beta \in \Lambda_{\alpha}, g\in \Gamma}[L_{\beta,g}, L_{-\beta,-g}]\subseteq H_{0}$$
and
$$\sum_{\beta \in \Lambda_{\alpha};\atop g, g^{'}\in \Gamma, g+g^{'}\neq 0}[L_{\beta,g}, L_{-\beta,g^{'}}]\subseteq \oplus_{g \in \Gamma\setminus \{0\} }H_{g}.$$

We also define
 $$V_{\Lambda_{\alpha}}:= \oplus_{\beta \in \Lambda_{\alpha}}L_{\beta}= \oplus_{\beta \in \Lambda_{\alpha},g \in \Gamma}L_{\beta, g}.$$

 Finally, we denote by $L_{\Lambda_{\alpha}}$ the following graded  subspace of $L$,

 $$L_{\Lambda_{\alpha}}:=H_{\Lambda_{\alpha}}\oplus V_{\Lambda_{\alpha}}.$$

\bprop\label{a8107217771}
For any $\alpha \in \Lambda$, the linear subspace $L_{\Lambda_{\alpha}}$ is a subalgebra of $L$.
\eprop

\bpf First we have to check that $L_{\Lambda_{\alpha}}$ satisfies  $[L_{\Lambda_{\alpha}}, L_{\Lambda_{\alpha}}]\subset L_{\Lambda_{\alpha}}.$      Taking into account $H=L_{o}$, we have
 \begin{equation}\label{810227}
[ L_{\Lambda_{\alpha}}, L_{\Lambda_{\alpha}}]=[H_{\Lambda_{\alpha}}\oplus V_{\Lambda_{\alpha}},H_{\Lambda_{\alpha}}\oplus V_{\Lambda_{\alpha}}] \subset [H_{\Lambda_{\alpha}}, V_{\Lambda_{\alpha}}] + [V_{\Lambda_{\alpha}}, H_{\Lambda_{\alpha}}] +\Sigma_{\beta,\delta \in \Lambda_{\alpha}}[L_{\beta},L_{\delta}].
\end{equation}

\noindent Let us consider the first summand in (\ref{810227}). Note that $H_{\Lambda_{\alpha}}$ is the direct sum of
$$\sum_{\beta \in \Lambda_{\alpha}, g\in \Gamma}[L_{\beta,g}, L_{-\beta,-g}]\subseteq H_{0}\subseteq L_{o}$$
and
$$\sum_{\beta \in \Lambda_{\alpha};\atop g, g^{'}\in \Gamma, g+g^{'}\neq 0}[L_{\beta,g}, L_{-\beta,g^{'}}]\subseteq \oplus_{g \in \Gamma\setminus \{0\} }H_{g}\subseteq L_{o}.$$

\noindent Given $\beta \in \Lambda_{\alpha}$, we have $[H_{\Lambda_{\alpha}}, L_{\beta}]\subset  L_{\beta \phi^{-1}}$, being $\beta \phi^{-1} \in \Lambda_{\alpha}$ by Lemma \ref{3556000}-2. Hence,
 \begin{equation}\label{8102275}
[H_{\Lambda_{\alpha}}, V_{\Lambda_{\alpha}}]\subset V_{\Lambda_{\alpha}}.
 \end{equation}

\noindent  Similarly, we can also get
 \begin{equation}\label{810227559}
 [V_{\Lambda_{\alpha}}, H_{\Lambda_{\alpha}}]\subset V_{\Lambda_{\alpha}}.
 \end{equation}

\noindent Consider now the third summand $\Sigma_{\beta,\delta \in \Lambda_{\alpha}}[L_{\beta},L_{\delta}]$.  Given $\beta, \delta \in \Lambda_{\alpha}$ such that $[L_{\beta},L_{\delta}]\neq 0$, if $\delta=-\beta$, then clearly $[L_{\beta},L_{\delta}]=[L_{\beta},L_{-\beta}]\subset H_{\Lambda_{\alpha}}.$ Suppose that $\delta\neq -\beta$. Since  $[L_{\beta},L_{\delta}]\neq 0$ together with  Lemma \ref{3556000}-2 ensures that  $\beta\phi^{-1}+\delta\phi^{-1} \in \Lambda$, we we have that $\{\beta,\delta\}$ is a connection from $\beta$ to $\beta\phi^{-1}+\delta\phi^{-1}$. The transitivity of $\sim$ gives now that $\beta\phi^{-1}+\delta\phi^{-1} \in \Lambda_{\alpha}$ and so
\begin{equation}\label{81022755958}
[L_{\beta},L_{\delta}]\subset L_{\beta\phi^{-1}+\delta\phi^{-1}}\subseteq V_{ \Lambda_{\alpha}}.
 \end{equation}

\noindent From Eqs. (\ref{810227})-(\ref{81022755958}) we conclude that $[L_{\Lambda_{\alpha}}, L_{\Lambda_{\alpha}}]\subset L_{\Lambda_{\alpha}}.$

Second, we have to verify that $\phi(L_{\Lambda_{\alpha}})= L_{\Lambda_{\alpha}}.$
But this is a direct consequence of Lemma \ref{3556000}-1 and Lemma \ref{3556jjsd}.
\epf

\bprop\label{a8107217771AFG}
 If $\gamma \not \in \Lambda_{\alpha}$ then $[L_{\Lambda_{\alpha}}, L_{\Lambda_{\gamma}}]=0$.

\eprop
\bpf
We have
\begin{equation}\label{810227acd}
[ L_{\Lambda_{\alpha}}, L_{\Lambda_{\gamma}}]=[H_{\Lambda_{\alpha}}\oplus V_{\Lambda_{\alpha}},H_{\Lambda_{\gamma}}\oplus V_{\Lambda_{\gamma}}] \subset [H_{\Lambda_{\alpha}}, V_{\Lambda_{\gamma}}] + [V_{\Lambda_{\alpha}}, H_{\Lambda_{\gamma}}] +[V_{\Lambda_{\alpha}},V_{\Lambda_{\gamma}}].
\end{equation}

\noindent Consider the above  third  summand $[V_{\Lambda_{\alpha}},V_{\Lambda_{\gamma}}]$ and suppose that there exist $\beta \in \Lambda_{\alpha}$ and $\eta \in \Lambda_{\gamma}$ such that $[L_{\beta}, L_{\eta}]\neq 0$. As necessarily $\beta \neq -\eta$, then $\beta\phi^{-1}+\eta\phi^{-1} \in \Lambda$. So $\{\beta, \eta, -\beta\phi^{-1}\}$ is a connection between $\beta$ and $\eta$. By the transitivity of the connection relation we have $\gamma \in \Lambda_{\alpha}$, a contradiction. Hence $[L_{\beta}, L_{\eta}]=0$ and so
\begin{equation}\label{caoyan}
[V_{\Lambda_{\alpha}},V_{\Lambda_{\gamma}}]=0.
\end{equation}

\noindent Consider now the first summand $[H_{\Lambda_{\alpha}}, V_{\Lambda_{\gamma}}]$ in (\ref{810227acd}) and suppose there exist
$\beta \in \Lambda_{\alpha}$ and $\eta \in \Lambda_{\gamma}$ such that $[[L_{\beta},L_{-\beta}], \phi(L_{\eta})]\neq 0.$ Then
$$[[L_{\beta,g},L_{-\beta,g^{'}}], \phi(L_{\eta})]\neq 0$$
for some $g, g^{'} \in \Gamma$. By Hom Jacobi identity, either $[L_{-\beta,g^{'}},\phi(L_{\eta})]\neq 0$ or $[L_{\beta,g},\phi(L_{\eta})]\neq 0$ and so $[V_{\Lambda_{\alpha}},V_{\Lambda_{\gamma}}]\neq 0$ in any case, what contradicts Eq. (\ref{caoyan}). Hence
$$ [H_{\Lambda_{\alpha}}, V_{\Lambda_{\gamma}}]=0.$$
Finally, we note that the same above argument shows
$$ [V_{\Lambda_{\gamma}}, H_{\Lambda_{\alpha}}]=0.$$
By Eq. (\ref{810227acd}) we conclude $[L_{\Lambda_{\alpha}}, L_{\Lambda_{\gamma}}]=0$.
\epf

\bthm\label{cao}
The following assertions hold.

$\rm 1$. For any $\alpha \in \Lambda$, the Hom-Lie  color subalgebra
 $$L_{\Lambda_{\alpha}}=H_{\Lambda_{\alpha}}\oplus V_{\Lambda_{\alpha}}$$
of $L$ associated to $\Lambda_{\alpha}$ is an ideal of $L$.

$\rm 2$. If $L$ is simple, then there exists a connection from $\alpha$ to $\beta$ for any $\alpha,\beta \in \Lambda$
and $H=\sum_{\alpha \in \Lambda}[L_{\alpha}, L_{-\alpha}]$.
\ethm

\bpf
$\rm 1$. Since $[L_{\Lambda_{\alpha}}, H]=[L_{\Lambda_{\alpha}}, L_{o}]\subset V_{[\alpha]}$,
taking into account Propositions \ref{a8107217771} and  \ref{a8107217771AFG}, we have
$$[L_{\Lambda_{\alpha}}, L]=[L_{\Lambda_{\alpha}}, H\oplus(\oplus_{\beta \in \Lambda_{\alpha}}L_{\beta})\oplus(\oplus_{\gamma \not \in \Lambda_{\alpha}}L_{\gamma})]\subset L_{\Lambda_{\alpha}}.$$
As we also have by  Proposition \ref{a8107217771} that $\phi(L_{\Lambda_{\alpha}} )=L_{\Lambda_{\alpha}} $, we conclude that $ L_{\Lambda_{\alpha}}$ is an ideal of $L$.

$\rm 2$. The simplicity of $L$ implies $L_{\Lambda_{\alpha}}=L$. From here, it is clear that  $\Lambda_{\alpha}=\Lambda$ and $H=\sum_{\alpha \in \Lambda}[L_{\alpha}, L_{-\alpha}]$.
\epf

\bthm\label{2hijianaas}
For a vector space complement $U$ of $span_{\mathbb{K}}\{{[L_{\alpha},L_{-\alpha}]: \alpha \in \Lambda}\}$ in H, we have
$$L = U + \sum\limits_{[\alpha] \in \Lambda/\sim}I_{[\alpha]},$$
where any $I_{[\alpha]}$ is one of the ideals  of $L$ described in Theorem \ref{cao}-1, satisfying
$[I_{[\alpha]},I_{[\beta]}]=0,$ whenever $[\alpha] \neq [\beta].$
\ethm

\bpf By Proposition \ref{alp93556jjsd77777}, we can consider the quotient set $\Lambda/\sim :=\{[\alpha]: \alpha \in \Lambda\}$. let us denote by
$ I_{[\alpha]}:= L_{\Lambda_{\alpha}}$. We have $I_{[\alpha]}$ is well defined and, by Theorem \ref{cao}-1, an ideal of $L$. Therefore
$$L = U + \sum\limits_{[\alpha] \in \Lambda/\sim}I_{[\alpha]}.$$
\noindent By applying Proposition \ref{a8107217771AFG} we also obtain $[I_{[\alpha]},I_{[\beta]}]=0$ if $[\alpha] \neq [\beta].$
\epf

\bdefn\label{2hijian}
The \textbf{annihilator} of a Hom-Lie color algebra $L$ is the set $\mathrm{Z}(L) = \{x \in L: [x, L]  = 0\}$.
\edefn

\bcor\label{667}
If $\mathrm{Z}(L) = 0$ and $[L,L] = L$, then $L$ is the direct sum of the ideals given in Theorem \ref{cao},
$$L =\oplus_{[\alpha] \in \Lambda/\sim}I_{[\alpha]}.$$
\ecor
\bpf
From $[L,L]=L$, it is clear that $L=\oplus_{[\alpha] \in \Lambda/\sim}I_{[\alpha]}.$ Finally, the sum is direct because $\mathrm{Z}(L) = 0$ and $[I_{[\alpha]}, I_{[\beta]}]=0$ if
 $[\alpha] \neq [\beta]$.
\epf

\section{The simplicity of split regular Hom-Lie color algebras of maximal length. }

In this section we focus on the simplicity of split regular Hom-Lie color algebras by centering our attention in those of maximal length. From now on char($\mathbb{K}$)=0.

\blem \label{lemma 4.3}
Let $L = H\oplus(\oplus_{\alpha \in \Lambda} L_{\alpha})$ be a split regular Hom-Lie color algebra. If $I$ is an ideal of $L$ then $I=(I\cap H)\oplus(\oplus_{\alpha \in \Lambda}(I\cap L_{\alpha})).$
\elem

\bpf
We may view $L = H\oplus(\oplus_{\alpha \in \Lambda} L_{\alpha})$ as a weight  module respect to  the split Hom-Lie color algebra $L_{0}$ with maximal abelian subalgebra $H_{0}$, (see Lemma \ref{1.2}-3), in  the natural way. The characteristic property of ideals gives us that $I$ is a submodule of $L$. It is well-known that a submodule of a weight module is again a weight module. From here, $I$ is a weight module respect to $L_{0}$, (and $H_{0}$), and so $I=(I\cap H)\oplus(\oplus_{\alpha \in \Lambda}(I\cap L_{\alpha})).$
\epf

Taking into account the above lemma, observe that the grading of $I$ and Lemma \ref{1.2}-1 let us write
\begin{equation}\label{zhangjian}
I=\oplus_{g \in \Gamma}I_{g}=\oplus_{g \in \Gamma}\big((I_{g}\cap H_{g})\oplus (\oplus_{\alpha \in \Lambda}(I_{g}\cap L_{\alpha,g})\big).
\end{equation}

\blem \label{lemma 4.1}
Let $L$ be a split regular Hom-Lie color algebra with $\mathrm{Z}(L)=0$ and $I$ an ideal of $L$.  If $I\subset H$ then $I=\{0\}$.
\elem
\bpf
Suppose there exists a nonzero ideal $I$ of $L$ such that $I \subset H$. We get $[I, H] \subset [H, H]=0$.  We also get $[I, \oplus_{\alpha \in \Lambda}L_{\alpha}]\subset I\subset H$. Then taking into account $H = L_{o}$, we have  $[I, \oplus_{\alpha \in \Lambda}L_{\alpha}]\subset H \cap (\oplus_{\alpha \in \Lambda}L_{\alpha})=0.$ From here $I\subset  \mathrm{Z}(L)=0$, which is a contradiction.
\epf

Let us  introduce the concepts of root-multiplicativity and maximal length in the framework of split Hom-Lie color algebras, in a similar way to the ones for split Hom Lie algebras (see \cite{BL52567}). For each $g \in \Gamma$, we denote by $\Lambda_{g}:=\{\alpha \in \Lambda: L_{\alpha,g}\neq 0\}.$

\bdefn
We say that a split regular Hom-Lie color algebra $L$ is \textbf{root-multiplicative} if given $\alpha \in \Lambda_{g_{i}}$ and $\beta \in \Lambda_{g_{j}}$, with $g_{i},g_{j}\in \Gamma$, such that $\alpha+\beta \in \Lambda$, then $[L_{\alpha,g_{i}},L_{\beta,g_{j}}]\neq 0$.
\edefn

\bdefn
We say that a split regular Hom-Lie color algebra $L$ is of \textbf{maximal length} if for any $\alpha \in \Lambda_{g}, g\in \Gamma$, we have $\mathrm{dim}L_{\kappa\alpha,\kappa g}=1$ for $\kappa \in \{\pm1\}$.
\edefn

Observe that if $L$ is of maximal lenth, then Eq. (\ref{zhangjian}) let us assert that given any nonzero ideal $I$ of $L$ then
\begin{equation}\label{zhangjian1}
I=\oplus_{g \in \Gamma}\big((I_{g}\cap H_{g})\oplus (\oplus_{\alpha \in \Lambda_{g}^{I}} L_{\alpha,g})\big).
\end{equation}
\noindent where $\Lambda_{g}^{I}:=\{\alpha \in \Lambda: I_{g}\cap L_{\alpha,g}\neq 0\}$ for each $g \in \Gamma$.

\bthm\label{THEOTEM 4.1 qw}
Let $L$ be a split regular Hom-Lie color algebra of maximal length, root multiplicative and with $\mathrm{Z}(L)=0$. Then $L$ is  simple if and only if it has all its nonzero roots connected and $H=\sum_{\alpha \in \Lambda}[L_{\alpha},L_{-\alpha}]$.
\ethm

\bpf
The first implication is Theorem \ref{cao}-2. To prove the converse,consider $I$ a nonzero ideal of $L$. By Lemma \ref{lemma 4.1} and Eq. (\ref{zhangjian1}) we can write $I=\oplus_{g \in \Gamma}\big((I_{g}\cap H_{g})\oplus (\oplus_{\alpha \in \Lambda_{g}^{I}} L_{\alpha,g})\big)$ with $\Lambda_{g}^{I} \subset \Lambda_{g}$ for any $g \in \Gamma$ and some $\Lambda_{g}^{I} \neq \emptyset$. Hence, we may choose $\alpha_{0} \in \Lambda_{g}^{I}$ being so
\begin{equation}\label{zhangjian12}
0\neq L_{\alpha_{0},g}\subset I.
\end{equation}
\noindent The fact $\phi(I)=I$ together with Lemma \ref{3556000}-1 allows us to assert that
\begin{equation}\label{zhangjian123}
  \mathrm{if} \quad  \alpha \in \Lambda_{I} \quad  \mathrm{then} \quad  \{\alpha\phi^{z}: z \in \mathbb{Z}\}\subset \Lambda_{I},
\end{equation}
\noindent that is
\begin{equation}\label{zhangjian1234}
\{L_{\alpha_{0}\phi^{z},g}: z \in \mathbb{Z}\} \subset I.
\end{equation}
\noindent Now, let us take any $\beta \in \Lambda$ satisfying $\beta \not \in \{\pm \alpha_{0}\phi^{z}: z \in \mathbb{Z}\}$. Since $\alpha_{0}$ and $\beta$ are connected, we have a connection $\{\gamma_{1},\gamma_{2},\cdots,\gamma_{k}\}$, $k\geq 2$, from $\alpha_{0}$ to $\beta$ satisfying:

 $\gamma_{1}= a_{0}\phi^{-n}$ for some $ n\in \mathbb{N}$, and

 $\gamma_{1}\phi^{-1}+\gamma_{2}\phi^{-1} \in \Lambda$,

  $\gamma_{1}\phi^{-2}+\gamma_{2}\phi^{-2}+\gamma_{3}\phi^{-1} \in \Lambda$,

 $\cdots \cdots \cdots$

 $\gamma_{1}\phi^{-i}+\gamma_{2}\phi^{-i}+\gamma_{3}\phi^{-i+1}+\cdots+\gamma_{i+1}\phi^{-1}\in \Lambda$,

 $\cdots \cdots \cdots$

 $\gamma_{1}\phi^{-k+2}+\gamma_{2}\phi^{-k+2}+\gamma_{3}\phi^{-k+3}+\cdots+\gamma_{i}\phi^{-k+i}+\cdots+\gamma_{k-1}\phi^{-1}\in \Lambda$,

 $\gamma_{1}\phi^{-k+1}+\gamma_{2}\phi^{-k+1}+\gamma_{3}\phi^{-k+2}+\cdots+\gamma_{i}\phi^{-k+i-1}+\cdots+\gamma_{k}\phi^{-1}= \epsilon \beta \phi^{-m}$ for some  $m \in \mathbb{N}$ and $\epsilon \in \{\pm 1\}.$

\noindent Consider $\gamma_{1}, \gamma_{2}$ and $\gamma_{1}+\gamma_{2}$. Since $\gamma_{2} \in \Lambda$ there exists $g_{1} \in \Gamma$ such that $L_{\gamma_{2}, g_{1}}\neq 0$. From here, the root-multiplicativity and maximal length of $L$ show $0 \neq [L_{\gamma_{1},g}, L_{\gamma_{2},g_{1}}]=L_{(\gamma_{1}+\gamma_{2})\phi^{-1},g+g_{1}}$, and by  Eq. (\ref{zhangjian1234})
$$ 0\neq L_{(\gamma_{1}+\gamma_{2})\phi^{-1},g+g_{1}}\subset I.$$
\noindent We can argue in a similar way from $\gamma_{1}\phi^{-1}+\gamma_{2}\phi^{-1}$, $\gamma_{3}$ and $\gamma_{1}\phi^{-2}+\gamma_{2}\phi^{-2}
+\gamma_{3}\phi^{-1}$ to get
$$0\neq L_{\gamma_{1}\phi^{-2}+\gamma_{2}\phi^{-2}
+\gamma_{3}\phi^{-1},g_{2}}\subset I$$
\noindent for some $g_{2} \in \Gamma$. Following this process with  the connection $\{\gamma_{1},\cdots,\gamma_{k}\}$ we obtain that
$$ 0\neq L_{\gamma_{1}\phi^{-k+1}+\gamma_{2}\phi^{-k+1}+\gamma_{3}\phi^{-k+2}+\cdots+\gamma_{k}\phi^{-1},g_{3}} \subset I$$
\noindent and so either $0\neq L_{\beta \phi^{-m},g_{3}}\subset I$ or $0\neq L_{-\beta \phi^{-m},g_{3}}\subset I$ for some $g_{3} \in \Gamma$. That is,
\begin{equation}\label{zhangjian1234888}
0\neq L_{\epsilon\beta \phi^{-m},g_{3}}\subset I \quad  \mathrm{for} \quad \mathrm{some} \quad \epsilon \in\{\pm 1\}, \quad \mathrm{some} \quad  g_{3}\in \Gamma
\end{equation}
\noindent and for any $\beta \in \Lambda$. By Lemma \ref{3556000}-1, we can get
\begin{equation}\label{zhangjian1234888}
0\neq L_{\epsilon\beta ,g_{3}}\subset I \quad  \mathrm{for} \quad \mathrm{some} \quad \epsilon \in\{\pm 1\}, \quad \mathrm{some} \quad  g_{3}\in \Gamma
\end{equation}
\noindent and for any $\beta \in \Lambda$.

 Taking into  account $H=\sum_{\gamma \in \Lambda}[L_{\gamma}, L_{-\gamma}]$, the grading of $L$ gives us $$H_{0}=\sum_{\gamma \in \Lambda, g \in \Gamma}[L_{\gamma,g},L_{-\gamma,-g}].$$ From here, there exists $\gamma \in \Lambda$ and $g_{4} \in \Gamma$ such that
\begin{equation}\label{zhangjian1234888600}
[[L_{\gamma,g_{4}},L_{-\gamma,-g_{4}}], \phi(L_{\epsilon\beta,g_{3}})]\neq 0.
\end{equation}
\noindent By the Hom Jacobi identity either $[L_{\gamma,g_{4}},\phi(L_{\epsilon\beta,g_{3}})]\neq 0$ or $[L_{-\gamma,-g_{4}},\phi(L_{\epsilon\beta,g_{3}})]\neq 0$ and so $L_{\gamma\phi^{-1}+\epsilon\beta\phi^{-2},g_{4}+g_{3}}\neq 0$ or
$L_{-\gamma\phi^{-1}+\epsilon\beta\phi^{-2},-g_{4}+g_{3}}\neq 0$. That is
\begin{equation}\label{zhangjian1234888600135}
0\neq L_{\kappa\gamma\phi^{-1}+\epsilon\beta\phi^{-2},\kappa g_{4}+g_{3}}\subset I
\end{equation}
\noindent for some $\kappa \in \{\pm 1\}.$ Since $\epsilon\beta  \in \Lambda_{g_{3}}$ we have by the maximal length of $L$ that $-\epsilon\beta  \in \Lambda_{-g_{3}}$. By Eq. (\ref{zhangjian1234888600135}) and the root-multiplicativity and maximal length of $L$ we obtain
\begin{equation}\label{zhangjian1234888600135qwe}
0\neq [ L_{\kappa\gamma\phi^{-1}+\epsilon\beta\phi^{-2},\kappa g_{4}+g_{3}}, L_{-\epsilon\beta\phi^{-2},-g_{3}}]=L_{\kappa\gamma\phi^{-2},\kappa g_{4}}\subset I.
\end{equation}
\noindent By Lemma \ref{3556000}-1, we can get
\begin{equation}\label{zhangjian1234888600135qwehky}
L_{\kappa\gamma,\kappa g_{4}}\subset I.
\end{equation}

\noindent Taking into account Eq. (\ref{zhangjian1234888600135qwe}) and that Eq. (\ref{zhangjian1234888600}) gives us
$$\beta\phi^{-1}([L_{\gamma,g_{4}},L_{-\gamma,-g_{4}}])\neq 0,$$
\noindent we have that for any $g_{5} \in \Gamma$ such that $L_{\epsilon\beta,g_{5}}\neq 0$ necessarily
$$0\neq[[L_{\gamma,g_{4}},L_{-\gamma,-g_{4}}], \phi(L_{\epsilon\beta,g_{5}})]=L_{\epsilon\beta\phi^{-1},g_{5}}\subset I$$
\noindent and so $L_{\epsilon\beta\phi^{-1}}\subset I$. That is, we can assert that
\begin{equation}\label{zh00135qwehky}
L_{\epsilon\beta}\subset I
\end{equation}

\noindent for any $\beta \in \Lambda$ and some $\epsilon \in \{\pm 1\}$. Since $H=\sum_{\beta \in \Lambda}[L_{\beta},L_{-\beta}]$  we get
\begin{equation}\label{zhawehky}
H\subset I.
\end{equation}

\noindent Now, given any $-\epsilon\beta \in \Lambda$, by the facts $-\epsilon\beta \neq 0$, $H\subset I$ and the maximal length of $L$ we have
\begin{equation}\label{zhawehky111}
[H_{0},L_{-\epsilon\beta}]=L_{-\epsilon\beta} \subset I.
\end{equation}

\noindent From Eqs. (\ref{zh00135qwehky})-(\ref{zhawehky111}) we conclude that $I=L$. Consequently $L$ is simple.
\epf

\bthm\label{THEOTEM 4.1}
Let $L$ be a split regular Hom-Lie color algebra of maximal length, root multiplicative and satisfying $\mathrm{Z}(L)=0$, $[L,L]=L$. Then $L$ is  the direct sum of the family of its minimal ideals, each one being a simple split regular Hom-Lie color algebra having all its nonzero roots connected.
\ethm
\bpf
By corollary \ref{667}, $L =\oplus_{[\alpha] \in \Lambda/\sim}I_{[\alpha]}$ is the direct sum of the ideals $I_{[\alpha]}=H_{\Lambda_{\alpha}}\oplus V_{\Lambda_{\alpha}}$=$(\sum_{\beta \in [\alpha]}[L_{\beta},L_{-\beta}])\oplus(\oplus_{\beta \in [\alpha]}L_{\beta})$ having  any $I_{[\alpha]}$ its root system, $\Lambda_{\alpha}$, with all of its roots connected. It is easy to check that $\Lambda_{\alpha}$ has
all of its roots $\Lambda_{\alpha}$-connected, (connected through roots in $\Lambda_{\alpha}$). We also have that any of the $I_{[\alpha]}$ is root-multiplicative as consequence of the root-multiplicativity  of $L$. Clearly $I_{[\alpha]}$ is of maximal length, and finally $Z_{I_{[\alpha]}}(I_{[\alpha]})$=0, (where  $Z_{I_{[\alpha]}}(I_{[\alpha]})$) denotes the center of $I_{[\alpha]}$ in  $I_{[\alpha]}$), as consequence of $[I_{[\alpha]}, I_{[\beta]}]=0 $ if $[\alpha]\neq [\beta]$, (Theorem \ref{2hijianaas}), and $Z(L)=0$. We can apply Theorem
\ref{THEOTEM 4.1 qw} to any $I_{[\alpha]}$  so as to conclude $I_{[\alpha]}$ is ssimple. It is clear that the decomposition $L =\oplus_{[\alpha] \in \Lambda/\sim}I_{[\alpha]}$ satisfies the assertions of the theorem.
\epf

\end{document}